\documentclass[a4paper,12pt]{amsart}
\usepackage{amsfonts}
\usepackage{amssymb}
\usepackage{ifthen}
\usepackage{graphicx}
\usepackage{color}
\nonstopmode \numberwithin{equation}{section}
\newtheorem{thm}{Theorem}
\newtheorem{lem}{Lemma}
\newtheorem{cor}{Corollary}
\newtheorem{prop}{Proposition}

\newtheorem{conj}{Conjecture}

\theoremstyle{definition}
\newtheorem{defn}{Definition}
\newtheorem{example}{Example}

\newtheorem{ques}{Question}
\newtheorem{rem}{Remark}

\newtheorem{rems}{Remarks}

\newcounter {own}
\def\theown {\thesection  .\arabic{own}}

\newenvironment{pf}[1][]{%
 \vskip 3mm
 \noindent
 \ifthenelse{\equal{#1}{}}%
  {{\slshape Proof. }}%
  {{\slshape #1.} }%
 }%
{\qed\bigskip}

\newcounter{alphabet}

\newcommand{\N}{{\mathbb N}}
\newcommand{\C}{{\mathbb C}}

\def\be{\begin{equation}}
\def\ee{\end{equation}}

\newcommand{\bee}{\begin{enumerate}}
\newcommand{\eee}{\end{enumerate}}

\newcommand{\blem}{\begin{lem}}
\newcommand{\elem}{\end{lem}}
\newcommand{\bthm}{\begin{thm}}
\newcommand{\ethm}{\end{thm}}
\newcommand{\bcor}{\begin{cor}}
\newcommand{\ecor}{\end{cor}}
\newcommand{\beg}{\begin{example}}
\newcommand{\eeg}{\end{example}}
\newcommand{\begs}{\begin{examples}}
\newcommand{\eegs}{\end{examples}}
\newcommand{\bdefn}{\begin{defn}}
\newcommand{\edefn}{\end{defn}}
\newcommand{\bprob}{\begin{prob}}
\newcommand{\eprob}{\end{prob}}
\newcommand{\bei}{\begin{itemize}}
\newcommand{\eei}{\end{itemize}}
\newcommand{\bqn}{\begin{ques}}
\newcommand{\eqn}{\end{ques}}
\newcommand{\bcon}{\begin{conj}}
\newcommand{\econ}{\end{conj}}
\newcommand{\bcons}{\begin{conjs}}
\newcommand{\econs}{\end{conjs}}
\newcommand{\bprop}{\begin{prop}}
\newcommand{\eprop}{\end{prop}}
\newcommand{\brem}{\begin{rem}}
\newcommand{\erem}{\end{rem}}
\newcommand{\brems}{\begin{rems}}
\newcommand{\erems}{\end{rems}}
\newcommand{\bo}{\begin{obser}}
\newcommand{\eo}{\end{obser}}
\newcommand{\bos}{\begin{obsers}}
\newcommand{\eos}{\end{obsers}}
\newcommand{\bpf}{\begin{pf}}
\newcommand{\epf}{\end{pf}}
\newcommand{\ba}{\begin{array}}
\newcommand{\ea}{\end{array}}
\newcommand{\beq}{\begin{eqnarray}}
\newcommand{\eeq}{\end{eqnarray}}
\newcommand{\beqq}{\begin{eqnarray*}}
\newcommand{\eeqq}{\end{eqnarray*}}

\newcommand{\ds}{\displaystyle}

\newcounter{minutes}
\divide\time by 60
\newcounter{hours}
\multiply\time by 60 \addtocounter{minutes}{-\time}

\begin{document}
\bibliographystyle{amsplain}
\title[Weighted composition operators between weighted Hardy spaces on rooted trees]
{Weighted composition operators between weighted Hardy spaces on rooted trees}


\author{P. Muthukumar}
\address{P. Muthukumar, Indian Statistical Institute,
Statistics and Mathematics Unit, 8th Mile, Mysore Road,
Bangalore, 560 059, India.}
\email{pmuthumaths@gmail.com}

\author{Ajay K.~Sharma}
\address{Ajay K.~Sharma, Department of Mathematics, Central University of Jammu,
Bagla, Rahya-Suchani, Samba 181143,
India.}
\email{aksju\_76@yahoo.com}

\author{Vivek Kumar}
\address{Vivek Kumar, Department of Mathematics, Central University of Jammu,
Bagla, Rahya-Suchani, Samba 181143,
India.}
\email{vivakkumar072@gmail.com}



\begin{abstract}
In this paper, we introduce a discrete analogue of weighted Hardy spaces
on rooted trees and study
weighted composition operators between them in detail. In particular, we
 characterize  bounded and compact  weighted composition operators between discrete  Hardy spaces.  We also consider isometric weighted composition operators between these spaces.
\end{abstract}
\thanks{
File:~\jobname .tex,
          printed: \number\day-\number\month-\number\year,
          \thehours.\ifnum\theminutes<10{0}\fi\theminutes
}
\maketitle
\pagestyle{myheadings}
\markboth{P. Muthukumar,  Ajay K. Sharma and Vivek Kumar}{Weighted composition operators
between weighted Hardy spaces on rooted trees}

\section{Introduction}\label{sec:intr}
Let $X$ be a linear space consisting of complex-valued maps on a non-empty set $\Omega$.
For a given self map $\varphi$ of $\Omega$ and a complex-valued map $\psi$ on $\Omega$, the corresponding weighted
composition operator $W_{\psi,\varphi}$ is defined by
$$ W_{\psi,\varphi}(f)= \psi (f\circ \varphi) \mbox {~for~} f \in X.
$$
If $\psi\equiv 1$, then $W_{\psi,\varphi}$ reduced to a composition operator $C_{\varphi}\,
($defined by $C_{\varphi}(f)=f\circ \varphi)$ and if $\varphi(x)=x$ for all $x\in \Omega$, then $W_{\psi,\varphi}$
reduced to a multiplication operator $M_{\psi}\, ($defined by $M_{\psi}(f)=\psi f)$.

The study of weighted composition operators were begun by Banach in early nineteenth century. The classical Banach-Stone theorem says that the surjective isometries between the spaces of continuous functions on some intervals, are certain weighted composition operators (see \cite{Banach}).
When $p\neq 2$, the isometric isomorphisms of the Hardy spaces $H^p$ and Bergman spaces  $A^p$
 are also weighted composition operators \cite{Forelli,Bergman isometry}.

Further, weighted composition operators appeared in other branches of mathematics
 such as dynamical systems and evolution equations. For example, refer \cite{Dynamical} to see the connection between weighted composition operators and classification of dichotomies in certain dynamical systems. Weighted composition transformations are also considered on some
compact convex family of analytic functions on the unit disk of the complex plane in order to solve some extremal problems in geometric function theory. For more details, refer \cite{arxiv W-CO,WCO}.

In recent years, studying operators on various function spaces defined on trees
becomes an interesting topic of research. See
\cite{Colonna-MO-3,Colonna-MO-2,Allen-MO-7,Colonna-MO-1,Ajay1} and references therein.
In \cite{MP-Tp-spaces}, the authors introduced the discrete Hardy spaces on trees and considered
multiplication operators on them. Composition operators on discrete Hardy spaces were
studied in \cite{CO-Tp-spaces,CO-Tp-spaces2}. Recently, multiplication operators between these function spaces were investigated in \cite{MSha1}. As an extension of all these works, in this paper, we introduce a discrete analogue of weighted Hardy spaces
on trees and study weighted composition operators between them.

  The paper is organized as follows. We refer to the Sections \ref{sec:prelim} and \ref{sec:basic} for basic definitions and preliminaries about the discrete weighted Hardy spaces $\mathbb{T}_{\sigma,p}$ such as inclusion properties, growth estimate, etc. In Sections \ref{sec:bdd} and \ref{sec:cmpt}, we  characterize  bounded and compact weighted composition operators between $\mathbb{T}_{\sigma,p}$ spaces, respectively.
  We discuss about isometric weighted composition operators between $\mathbb{T}_{\sigma,p}$ spaces in Section \ref{sec:iso}.  Finally, in Section \ref{sec:eg}, we compare the boundedness and compactness of the operators $M_{\psi}, C_{\varphi}$ with the boundedness and compactness of the operator $W_{\psi,\varphi}$.

 \section{Preliminaries}\label{sec:prelim}

 Let $\mathbb{N}, \mathbb{N}_0$ and $\mathbb{C}$ denote the set of all positive
 integers, the set of all non-negative integers and the set of all complex numbers,
 respectively.  An infinite tree (graph) $T$ with a special vertex (known as a root),
 is called a \textit{rooted tree}. Using the edge counting distance, $T$ can be regarded
 as a metric space. For the basic graph theory, we refer the readers to \cite{Book:graph}.
 By abuse of language, we denote $x\in T$ when $x$ is in the
 vertex set of $T$ and  we call a function defined on the vertex set of $T$ by a
 function defined on $T$. For each vertex $x\in T$, $|x|$ denotes the distance between
 the root $\textsl{o}$ and the vertex $x$. For $n\in \N_0$, let $D_n$ and $c_n$ denote the
 set of all vertices $x\in T$ with $|x|=n$ and the number  of vertices $x\in T$ with
 $|x|=n$, respectively.

A positive-valued function $T$  is called a \textit{weight function} on $T$. Fix a weight
function $\sigma$ on $T$ and fix $p>0$. For $n\in \N_0$ and a complex-valued function $f$ on
$T$, we introduce
$$
M_{\sigma,p}(n,f)=\left(\frac{1}{c_n} \sum\limits_{|x|=n} \sigma^p(x)|f(x)|^p \right)^{1/p}.
$$
Now, we define discrete weighted Hardy spaces $\mathbb{T}_{\sigma, p}$ as follows.
For $p>0$,
$$
\mathbb{T}_{\sigma,p}=\left\{f:T\rightarrow \C : \|f\|_{\sigma,p}:=\sup\limits_
{n\in \N_0} M_{\sigma,p}(n,f) <\infty \right\},
$$
and
$$
\mathbb{T}_{\sigma,\infty}=\left\{f:T\rightarrow \C : \|f\|_{\sigma,\infty}:=
\sup\limits_{x\in T} \sigma(x)|f(x)|<\infty \right\}.
$$

In particular, if $\sigma\equiv 1$, then the space $\mathbb{T}_{\sigma,p}$ becomes
the discrete Hardy space $\mathbb{T}_{p}$ with $\|\cdot\|_p= \|\cdot\|_{1,p}$ (see \cite{MP-Tp-spaces}). For a complex-valued
function $f$ on $T$, it is evident that $\|f\|_{\sigma,p}=\|\sigma f\|_{p}$. Thus,
$f\in \mathbb{T}_{\sigma,p}$ is equivalent to $\sigma f\in \mathbb{T}_p$.

Throughout the paper,  $T$ denotes an infinite rooted tree, $\sigma$ denotes a weight function on $T$.

\section{Basic Properties}\label{sec:basic}
 We begin with the inclusion properties of  $\mathbb{T}_{\sigma,p}$ spaces.
\bprop\label{inc}

\begin{enumerate}
\item If $\sup_{n\in \mathbb{N}_0} c_n= \infty$, then $\mathbb{T}_{\sigma,q}\subseteq\mathbb{T}_{\sigma,p}$ for $0<p<q$.
 \item If $\sup_{n\in \mathbb{N}_0} c_n< \infty$, then $\mathbb{T}_{\sigma,p}=\mathbb{T}_{\sigma,\infty}$ for all $p>0$.
\end{enumerate}
\eprop
Proof of this proposition is similar to the case of $\mathbb{T}_p$ spaces.
See \cite[Section 3]{MSha1} for a proof of later case. Also, see Corollary \ref{incl}
below for different weight functions.
In view of Proposition \ref{inc}, throughout this article,
we assume that $\{c_n: n\in \N_0\}$ is unbounded to avoid the triviality.

\bthm
For $1 \leq p \leq\infty, \mathbb{T}_{\sigma,p}$ is a Banach space.
\ethm
As the proof is straight forward, we omit the details here. Refer
\cite[Theorem 3.1]{MP-Tp-spaces} for case of $\mathbb{T}_p$ spaces.

Now, fix $x\in T$. For any $f\in \mathbb{T}_{\sigma,p}$, we obtain
$$
\sigma^p(x)|f(x)|^p\leq \sum\limits_{|y|=n} \sigma^p(y)|f(y)|^p\leq c_n \|f\|^p_{\sigma,p},
$$
where $n=|x|$. Thus, we get the following growth estimate.

\bprop\label{growth}
 For every $x\in T$, we have
$$
|f(x)|\leq\frac{\left(c_{|x|}\right)^{1/p}}{\sigma(x)}\|f\|_{\sigma,p}, \,\,  f\in \mathbb{T}_{\sigma,p}.
$$
Sharpness of this estimate can be easily verified by using characteristic functions.
\eprop

\bcor
For every $x\in T$, the evaluation map $ev_x:\mathbb{T}_{\sigma,p}\rightarrow\C$ defined
 by $ev_x(f):=f(x)$ is a bounded linear functional on $\mathbb{T}_{\sigma,p}$.
\ecor

\bcor
For $1 \leq p \leq\infty, \mathbb{T}_{\sigma,p}$ is a functional Banach space.
\ecor

Recall that Banach space $X$ consisting of functions defined on a set $\Omega$
is called \textit{functional Banach space} if for each $x\in \Omega$, the point
evaluation map $ev_x : f\in X \mapsto f(x)$ is a bounded linear functional on $X$
and $f(x)=f(y)$ for each $f\in X$ implies $x=y.$

\bthm
$\mathbb{T}_p \subseteq \mathbb{T}_{\sigma,p}$ if and only if $\sigma$ is a bounded map.
\ethm

\bpf
 Suppose $\sigma$ is a bounded map and $f\in \mathbb{T}_p.$ From the definition of
 $\|\cdot\|_p$, it is trivial to see that
 $$
 \|\sigma f\|_p\leq \|\sigma\|_\infty \|f\|_p.
 $$
This yields that $\sigma f\in \mathbb{T}_p,$ that is, $f\in \mathbb{T}_{\sigma,p}.$
Therefore, $\mathbb{T}_p \subseteq \mathbb{T}_{\sigma,p}.$

For the converse part, we suppose that $\sigma$ is unbounded.  For each $n\in \mathbb N,$
choose $v_n\in T$ such that $\sigma(v_n)>n.$ Define $f : T \to \mathbb C$ by
\[ f(v) = \left\{
 \begin{array}{cl}
     \left(c_{|v|}\right)^{1/p}&\mbox{ if }v=v_n \mbox{ for some } n\in \mathbb N,\\
0  & \mbox{ otherwise}.
\end{array}
\right.
\]
Thus, \[ M_p(k, f) = \left\{
 \begin{array}{ll}
   1  & \mbox{ if }  k = v_n  \mbox{ for some }  n\in \mathbb N,\\
0  & \mbox{ otherwise}.
\end{array}
\right.
\]
Therefore, $f\in \mathbb{T}_p$ with $\|f\|_p=1.$
It is easy to see that
\begin{align*}
\|f\|_{\sigma,p} & =\sup_{n\in \mathbb N} \frac{\sigma(v_n) |f(v_n)|}{\left(c_{|v_n|}\right)^{1/p}}  = \sup_{n\in \mathbb N} \sigma(v_n),
\end{align*}
which is not finite.  Hence $f \notin \mathbb{T}_{\sigma,p}.$
This completes the proof.
\epf

\bthm
 $\mathbb{T}_{\sigma,p} \subseteq \mathbb{T}_p$ if and only if $\sigma$ is bounded away
 from zero, that is,  there is a constant $m>0$ such that $m \leq \sigma(x)$   {for all}
  $ x\in T.$
\ethm

\bpf
Suppose that $\sigma$ is bounded away from zero. Then, $1/\sigma$ is a
 bounded function.
Let $f\in \mathbb{T}_{\sigma,p}$, i.e., $\sigma f\in \mathbb{T}_p.$ Consequently, $f=\left(1/\sigma\right)(\sigma f)\in \mathbb{T}_p$.
Hence $\mathbb{T}_{\sigma,p} \subseteq \mathbb{T}_p.$

For the other way inclusion, suppose that $\sigma$ is not away from zero. Choose
 $v_n\in T$ such that $\sigma(v_n)< 1/n$ for all $n\in \mathbb N.$
Define $f$ by
\[
 f(v) = \left\{
 \begin{array}{cl}
   n\left(c_{|v_n|}\right)^{1/p}   & \mbox{ if }  v=v_n  \mbox{ for some }  n\in \mathbb N,\\
0  & \mbox{ otherwise}.
\end{array}
\right.
\]
As $M_p(|v_n|,f)=n \to\infty$, we have $f \notin \mathbb{T}_p.$
But,
 \begin{align*}
 \|f\|_{\sigma,p}  =\sup_{n\in \mathbb N}M_{\sigma,p} (|v_n|,f)  = \sup_{n\in \mathbb N} \frac{\sigma(v_n)|f(v_n)|}{\left(c_{|v_n|}\right)^{1/p}} \leq 1,
 \end{align*}
that is, $f\in \mathbb{T}_{\sigma,p}$. The desired result follows.
\epf

\bcor
$\mathbb{T}_p=\mathbb{T}_{\sigma,p}$ if and only if there exist $m,M>0$ such that
$0<m\leq\sigma(x)\leq M \text{ for all } x\in T.$
In this case,
$$
m\|f\|_p \leq \|f\|_{\sigma,p}\leq M \|f\|_p \text{ for all } f\in \mathbb{T}_p   (=\mathbb{T}_{\sigma,p}).
$$
\ecor

\brem
The inclusion $\mathbb{T}_p \subseteq \mathbb{T}_{\sigma,p}$ is equivalent to saying that
$\sigma f\in \mathbb{T}_p$ whenever $f\in \mathbb{T}_p$; that is, the multiplication
operator $M_\sigma$ maps $\mathbb{T}_p$ to $\mathbb{T}_p$. In turn, this is same as
$\sigma$ is bounded (see \cite[Theorem D]{MSha1}). Similarly, the inclusion
$\mathbb{T}_{\sigma,p} \subseteq \mathbb{T}_p$ is equivalent to the multiplication
operator $M_{1/\sigma}$ maps $\mathbb{T}_p$ to $\mathbb{T}_p$.
\erem

\section{Bounded weighted composition operators}\label{sec:bdd}
 Throughout this article, $\sigma_1, \sigma_2,$ denote weight functions on $T$,
 $\varphi$ denotes a self map of $T$ and $\psi$ denotes a complex-valued map on $T$.
In this section, we will discuss about the boundedness of the weighted composition operator
$W_{\psi,\varphi}:f \mapsto \psi (f\circ\varphi)$ between $\mathbb{T}_{\sigma,p}$ spaces.
When we discuss an operator $A:X\to Y$, it is understood that the operator norm
$\|A\|_{X\to Y}$ will be simply denoted by $\|A\|$.

\bthm\label{infty-p}
For $0<p\leq \infty$, the operator  $W_{\psi,\varphi}:\mathbb{T}_{\sigma_1,\infty}\to \mathbb{T}_{\sigma_2,p}$ is bounded if and only if $\Psi\in \mathbb{T}_{\sigma_2,p}.$
In this case,
$$
\|W_{\psi,\varphi}\| =\left\|\Psi\right\|_{\sigma_2,p},
$$
where $\Psi=\psi/(\sigma_1\circ\varphi)$.
\ethm

\bpf
 Suppose that $\Psi\in \mathbb{T}_{\sigma_2,p}.$ First,
 we consider the case $p=\infty$.
 For $x\in T$ and $f\in \mathbb{T}_{\sigma_1,\infty},$
\begin{align*}
\sigma_2(x)|\psi(f\circ\varphi)(x)|
 =\sigma_2(x)  \frac{|\psi(x)|}{\sigma_1(\varphi(x))} \sigma_1(\varphi(x))
 |f(\varphi(x))|
 \leq \left\| \Psi\right\|_{\sigma_2,\infty}
 \|f\|_{\sigma_1,\infty}.
\end{align*}
Therefore,
$$
\|W_{\psi,\varphi}f\|_{\sigma_2,\infty} \leq \left\| \Psi\right\|_{\sigma_2,\infty} \|f\|_{\sigma_1,\infty} \text{ for all } f\in \mathbb{T}_{\sigma_1,\infty}.
$$

 Now, we consider the case $p<\infty$.
 For each $n\in \N_0$,
\begin{align*}
 M_{\sigma_2,p}^p(n,W_{\psi,\varphi}f)  =\frac{1}{c_n} \sum_{|x|=n} \frac{\sigma_2^p(x) |\psi(x)|^p}{\sigma_1^p(\varphi(x))} |f(\varphi(x))|^p \sigma_1^p(\varphi(x))
  \leq \|f\|_{\sigma_1,\infty}^p \left\| \Psi\right\|_
 {\sigma_2,p }^p.
 \end{align*}
Thus,
$$
\| W_{\psi,\varphi}f\|_{\sigma_2,p} \leq \left\| \Psi\right\|_{\sigma_2,p} \|f\|_{\sigma_1,\infty} \text{ for all } f\in \mathbb{T}_{\sigma_1,\infty}.
$$

Hence, $ W_{\psi,\varphi}:\mathbb{T}_{\sigma_1,\infty}\to \mathbb{T}_{\sigma_2,p}$ is a bounded operator with
$$\| W_{\psi,\varphi}\| \leq \left\| \Psi\right\|_{\sigma_2,p}.
$$

Conversely, suppose that $ W_{\psi,\varphi}:\mathbb{T}_{\sigma_1,\infty} \to \mathbb{T}_{\sigma_2,p}$ is a bounded operator.  Then,
$$
\Psi= W_{\psi,\varphi}\left(1/\sigma_1\right)\in \mathbb{T}_{\sigma_2,p}\;  \text{ as } \; 1/\sigma_1\in \mathbb{T}_{\sigma_1,\infty}.
$$
Since $\displaystyle\left\|1/\sigma_1\right\|_{\sigma_1,\infty}=1,$
 we also have that
$$
\| W_{\psi,\varphi}\|\geq \left\| \Psi\right\|_{\sigma_2,p}.
$$
Thus, the desired result follows.
\epf

Upon taking $\psi(x)=1$ and $\varphi(x)=x,\, x\in T$, in Theorem \ref{infty-p}, we get
the following results, respectively.

\bcor\label{cor}
\begin{itemize}
 \item[(i)] Let  $0<p\leq\infty$. Then, the composition operator $C_{\varphi}:\mathbb{T}_{\sigma_{1},\infty}\rightarrow \mathbb{T}_{\sigma_{2},p}$
      is bounded if and only if $1/(\sigma_{1}\circ\varphi)\in \mathbb{T}_{\sigma_{2},p}$.
    In this case,
       $$
       \|C_{\varphi}\|=\left\|1/(\sigma_1\circ \varphi)\right\|_{\sigma_{2},p}.
       $$
\item[(ii)] Let $0<p\leq\infty$. Then, the multiplication operator $M_{\psi}:\mathbb{T}_{\sigma_{1},\infty}\rightarrow \mathbb{T}_{\sigma_{2},p}$
    is bounded if and only if
    $\psi/\sigma_{1}\in \mathbb{T}_{\sigma_{2},p}.$
    Moreover,
    $$
    \|M_{\psi}\|=\left\|\psi/\sigma_{1}\right\|_{\sigma_{2}, p}.
    $$
\end{itemize}
\ecor

\bthm\label{p-infty}
 For $0<p<\infty$, the operator $W_{\psi,\varphi}:\mathbb{T}_{\sigma_1,p } \to \mathbb{T}_{\sigma_2,\infty}$ is  bounded if and only if
 $$
 M = \sup\limits_{x\in T}\left| \frac{\sigma_2\psi}{\sigma_1\circ\varphi}(x)\right| \left(c_{|\varphi(x)|}\right)^{1/p}<\infty.
 $$
 Furthermore, $\|W_{\psi,\varphi}\|= M.$
 \ethm

\bpf
Suppose that $ W_{\psi,\varphi} : \mathbb{T}_{\sigma_1,p}  \to \mathbb{T}_{\sigma_2,\infty}$ is a bounded operator. Fix $x\in T.$  Define $f$ by
\[ f(y) = \left\{
 \begin{array}{cl}
   \left(c_{|x|}\right)^{1/p}\ds/\sigma_1(x)   & \mbox{ if }  y=x,\\
0  & \mbox{ otherwise}.
\end{array}
\right.
\]
Then, $\|f\|_{\sigma_1,p}=1$ and
$$
\|W_{\psi,\varphi}\|\geq \|W_{\psi,\varphi}(f)\|_{\sigma_2,\infty} \geq \frac{\sigma_2(x)|\psi(x)|}{\sigma_1(\varphi(x))}  c_{|\varphi(x)|}^{1/p}.
$$
Since $x$ is arbitrary, it is evident that
$$
\|W_{\psi,\varphi}\|\geq \sup\limits_{x\in T}\left| \frac{\sigma_2\psi}
{\sigma_1\circ\varphi}(x)\right| \left(c_{|\varphi(x)|}\right)^{1/p}.
$$

Conversely, suppose that
$$
M=\sup_{x\in T}\frac{\sigma_2(x)|\psi(x)|}{\sigma_1(\varphi(x))}
 c_{|\varphi(x)|}^{1/p} <\infty.
 $$
For $f\in \mathbb{T}_{\sigma_1,p},$  by the growth estimate (Proposition \ref{growth}),
 we have
$$
\sigma_2(x)|\psi(x)f(\varphi(x))| \leq \sigma_2(x) |\psi(x)| \frac{c_{|\varphi(x)|}^{1/p}\|f\|_{\sigma_1,p}}{\sigma_1(\varphi(x))}.
$$
Thus,
$$
 \sigma_2(x) |W_{\psi,\varphi}f(x)| \leq M \,\|f\|_{\sigma_1,p} \; \text { for all }  x\in T.
 $$
Therefore, $\| W_{\psi,\varphi}f\|_{\sigma_2,\infty}\leq M\, \|f\|_{\sigma_1,p}$
  for all $f\in \mathbb{T}_{\sigma_1,p}.$
This yields that,
 $ W_{\psi,\varphi}: \mathbb{T}_{\sigma_1,p} \to \mathbb{T}_{\sigma_2,\infty}$
 is a bounded operator with $\| W_{\psi,\varphi}\| \leq M.$
This completes the proof.
\epf

As a special case of Theorem \ref{p-infty}, we have the following results.
\bcor
\begin{itemize}
  \item[(i)] For $0<p<\infty$, the composition operator
$C_{\varphi}:\mathbb{T}_{\sigma_{1},p}\rightarrow \mathbb{T}_{\sigma_{2},\infty}$
     is bounded if and only if
$$
\|C_{\varphi}\|=\sup\limits_{x\in T} \frac{\sigma_2(x)}{\sigma_1(\varphi(x))} \left(c_{|\varphi(x)|}\right)^{1/p}<\infty.
$$

\item[(ii)] For $0<p<\infty$, the multiplication operator
$M_{\psi}:\mathbb{T}_{\sigma_{1},p}\rightarrow \mathbb{T}_{\sigma_{2},\infty}$
 is bounded if and only if
$$
\|M_{\psi}\|=\sup\limits_{x\in T} \frac{\sigma_2(x)|\psi(x)|}{\sigma_1(x)} \left(c_{|x|}\right)^{1/p}<\infty.
$$

\end{itemize}

\ecor

Before we discuss about boundedness of the operator
$ W_{\psi,\varphi}: \mathbb{T}_{\sigma_1,p} \to \mathbb{T}_{\sigma_2,p}$,
let us introduce few notations. The characteristic (indicator) function at a vertex
$x\in T$ will be denoted by $\chi_{x}$.
Recall that
$D_n=\{x\in T: |x|=n\}, n\in \N_0$. We denote the set $\{|\varphi(x)|: x\in D_n\}$
by  $[\varphi(D_n)]$.
For $w\in T, n\in \N_0$, let $ N_{\varphi}(n,w)$ denote the number of elements in $\varphi^{-1}(w)\cap D_n$ (the number of elements in $D_n$ that are mapped to $w$
under $\varphi$). For  $m,n\in \N_0$, define
$$
N_{m,n}:= \displaystyle \max_{|w|=m} N_{\varphi}(n,w),
$$
 the maximum number of elements in $D_n$ that are mapped to a single element in $D_m.$
For $m,n\in \mathbb N_0,$ there exist $v_{m,n}\in D_m$ such that $N_\varphi(n,v_{m,n})=N_{m,n},$ equivalently, there are $N_{m,n}$ vertices of $D_n$ those are mapped into $v_{m,n}\in D_m$ under $\varphi.$

\bthm\label{p-p}
Let $0<p<\infty$ and for $m,n\in \mathbb N_0,$ choose $v_{m,n}\in D_m$ as above.
\begin{itemize}
  \item[(i)] If
$W_{\psi,\varphi}: \mathbb{T}_{\sigma_1,p} \to \mathbb{T}_{\sigma_2,p}$
 is a bounded operator, then
 $$
\sup\limits_{n\in \mathbb N_0}\left(\frac{1}{c_n} \sum_{m\in [\varphi(D_n)]} \sum\limits_{\substack{|x|=n\\\varphi(x)=v_{m,n}}} \left|\frac{\sigma_2 \psi}
{\sigma_1 \circ \varphi} (x)\right|^p c_m \right)\leq \|W_{\psi,\varphi}\|^p.
$$
  \item[(ii)] If
   $$
  M=\sup_{n\in \mathbb N_0}\left(\frac{1}{c_n} \sum_{m\in [\varphi(D_n)]}
  \sum_{\substack{|x|=n\\ |\varphi(x)|=m}} \left|\frac{\sigma_2
   \psi}{\sigma_1 \circ \varphi} (x)\right|^p c_m\right) <\infty,
   $$
   then $W_{\psi,\varphi}: \mathbb{T}_{\sigma_1,p} \to \mathbb{T}_{\sigma_2,p}$
   is a bounded operator with $\|W_{\psi,\varphi}\|^p \leq M$.
\end{itemize}
\ethm

\bpf Suppose $W_{\psi,\varphi}: \mathbb{T}_{\sigma_1,p} \to \mathbb{T}_{\sigma_2,p}$
 is a bounded operator.
 For each $n\in \mathbb N_0,$ consider a function $f_n$ on $T$, defined  by
 $$
 f_n(x)=\sum_{m=0}^{\infty} \frac{c_m^{1/p}}{\sigma_1(v_{m,n})} \chi_{v_{m,n}}(x).
 $$
Thus, $M_{\sigma_1,p}(k,f_n)=1$ for all $k, n\in \mathbb N_0$ and hence
$\|f_n\|_{\sigma_1,p}=1$ all $n\in \mathbb N_0$.
Note that $f_n(\varphi(x))\neq 0$ and $|\varphi(x)|=m$ for some $m\in \mathbb N_0$,
forces that $\varphi(x)=v_{m,n}$.
Now,
\begin{align*}
 M_{\sigma_2,p}^{p}(n,\psi(f_n \circ \varphi))
 & = \frac{1}{c_n} \sum_{m=0}^{\infty} \sum_{\substack{|x|=n\\|\varphi(x)|=m}} \sigma_2^p(x)|\psi(x)|^p |f_n(\varphi(x))|^p \\
& =\frac{1}{c_n}\sum_{m\in [\varphi(D_n)]} \sum_{\substack{|x|=n\\\varphi(x)=v_{m,n}}} \left|\frac{\sigma_2\psi}{\sigma_1 \circ \varphi} (x)\right|^p  c_m.
 \end{align*}
As
$M_{\sigma_2,p}^p (n,\psi(f_n\circ \varphi)) \leq \|\psi(f_n\circ\varphi)\|^p_
{\sigma_1,p}\leq \|W_{\psi,\varphi}\|^p \text{ for all } n\in \mathbb N_0,
$
we see that
$$
\sup\limits_{n\in \mathbb N_0}\left(\frac{1}{c_n} \sum_{m\in [\varphi(D_n)]} \sum\limits_{\substack{|x|=n\\\varphi(x)=v_{m,n}}} \left|\frac{\sigma_2 \psi}
{\sigma_1 \circ \varphi} (x)\right|^p c_m \right)\leq \|W_{\psi,\varphi}\|^p.
$$
Hence the result $\textrm{(i)}$ follows.

 Now, we prove the result $\textrm{(ii)}$. Fix $n\in \mathbb N_0.$ For any
 $f\in \mathbb{T}_{\sigma_1,p}$, growth estimate gives that
 \begin{align*}
  M_{\sigma_2,p}^p(n,\psi f \circ \varphi) & =\frac{1}{c_n} \sum_{m\in [\varphi(D_n)]} \sum_{\substack{ |x|=n\\ |\varphi(x)|=m}} \sigma_2^p(x)|\psi(x)|^p|f(\varphi(x))|^p \\
& \leq \frac{1}{c_n} \sum_{m\in [\varphi(D_n)]} \sum_{\substack{ |x|=n\\|\varphi(x)|=m}} \sigma_2^p(x)|\psi(x)|^p \frac{c_m\|f\|_{\sigma_1,p}^p}{\sigma_1^p(\varphi(x))}.
\end{align*}
Thus, for every  $f\in \mathbb{T}_{\sigma_1,p}$, one has
$$
\|W_{\psi,\varphi}(f)\|^p_{\sigma_2,p} \leq
\left( \sup_{n\in \mathbb N_0} \frac{1}{c_n} \sum_{m\in [\varphi(D_n)]}
 \sum_{\substack{ |x|=n\\|\varphi(x)|=m}} \left|\frac{\sigma_2 \psi}
{\sigma_1 \circ \varphi} (x)\right|^p c_m \right)\|f\|_{\sigma_1,p}^p .
$$
This completes the proof of $\textrm{(ii)}$.
\epf

\brem
Though we do not have a necessary and sufficient condition for the boundedness of the operator $W_{\psi,\varphi}: \mathbb{T}_{\sigma_1,p} \to \mathbb{T}_{\sigma_2,p}$, the upper and lower
bounds for $\|W_{\psi,\varphi}\|$ that we get in Theorem \ref{p-p}, are very close to each other.
\erem

\bprop
Let $0<p<\infty$. If $$
\sup_{n\in \mathbb N_0} \frac{1}{c_n} \sum_{m=0}^\infty N_{m,n} c_m \left(
 \sup_{\substack{|x|=n\\|\varphi(x)|=m}} \left|\frac{\sigma_2 \psi}{\sigma_1
\circ \varphi} (x)\right|^p \right) <\infty,
$$
 then $W_{\psi,\varphi}:\mathbb{T}_{\sigma_1,p} \to \mathbb{T}_{\sigma_2,p}$
 is a bounded operator.
\eprop

\bpf
Fix $n\in \mathbb N_0$ and $f\in \mathbb{T}_{\sigma_1,p}.$ Then,
\begin{align*}
  M_{\sigma_2,p}^p(n,\psi(f\circ\varphi)) & =\frac{1}{c_n} \sum_{m=0}^\infty \sum_{\substack{|x|=n\\|\varphi(x)|=m}} \sigma_2^p(x)|\psi(x)|^p |f(\varphi(x))|^p \frac{\sigma_1^p(\varphi(x))}{\sigma_1^p(\varphi(x))}\\
   & \leq \frac{1}{c_n} \sum_{m=0}^\infty \sup_{\substack{|x|=n\\|\varphi(x)|=m}} \left|\frac{\sigma_2\psi}{\sigma_1 \circ \varphi} (x)\right|^p \sum_{\substack
   {|x|=n\\|\varphi(x)|=m}} \sigma_1^p(\varphi(x))|f(\varphi(x))|^p.
  \end{align*}
Since every element of $D_m$ has atmost $N_{m,n}$ pre-images in $D_n,$ we have
$$
M_{\sigma_2,p}^p(n,\psi(f\circ\varphi)) \leq \frac{1}{c_n} \sum_{m=0}^\infty \sup_{\substack{|x|=n\\|\varphi(x)|=m}} \left|\frac{\sigma_2 \psi}{\sigma_1\circ\varphi} (x)\right|^p  N_{m,n}\, c_m\|f\|_{\sigma_1,p}^p.
$$
This yields that
$W_{\psi,\varphi}:\mathbb{T}_{\sigma_1,p} \to \mathbb{T}_{\sigma_2,p}$ is a bounded
operator with
$$ \|W_{\psi,\varphi}\|^p \leq \sup_{n\in \mathbb N_0} \left(\frac{1}{c_n} \sum_{m=0}^\infty \sup_{\substack{|x|=n\\|\varphi(x)|=m}} \left|\frac{\sigma_2\psi}{\sigma_1\circ\varphi}
 (x)\right|^p  N_{m,n}\, c_m\right).
$$
\epf

\brem
If we take $\sigma_1= \sigma_2= \psi \equiv 1$, then
$\mathbb{T}_{\sigma_1,p}=\mathbb{T}_{\sigma_2,p}=\mathbb{T}_p$ and
$W_{\psi,\varphi}=C_\varphi.$ As a special case of the previous theorem,
we get
$$
\|C_\varphi\|^p \leq \sup_{n\in \mathbb N_0}\frac{1}{c_n}
\sum_{m=0}^\infty N_{m,n}\,c_m.
$$
Indeed,
$$\|C_\varphi\|^p =\sup_{n\in \mathbb N_0}\frac{1}{c_n}
 \sum_{m=0}^\infty N_{m,n}\,c_m.
 $$
 See \cite[Theorem 4]{CO-Tp-spaces2} for a proof.
\erem

\bthm\label{bdd-Mo}
Let $0<p\leq\infty$. Then, the multiplication operator $M_{\psi}:\mathbb{T}_{\sigma_{1},p}\rightarrow \mathbb{T}_{\sigma_{2},p}$  is bounded
 if and only if
$\psi/\sigma_{1}\in \mathbb{T}_{\sigma_{2},\infty}$. Further, we also have
$$
\|M_{\psi}\|=\left\|\psi/\sigma_{1}\right\|_{\sigma_{2},\infty}.
$$
\ethm

\bpf
Suppose the operator
 $M_{\psi}:\mathbb{T}_{\sigma_{1},p}\rightarrow \mathbb{T}_{\sigma_{2},p}$
 is bounded. For the case $p=\infty$, see Corollary \ref{cor}. Now, assume that
 $p<\infty$.  Fix $y\in T.$
Define $$
f_{y}(x) = \left\{
 \begin{array}{cl}
   \frac{c_{|y|}^{ 1/p}}{\sigma_{1}(y)}  & \mbox{ if }  x=y,\\
0  & \mbox{ otherwise}.
\end{array}
\right.
$$
Then, it is trivial to see that $\|f_{y}\|_{\sigma_{1},p}=1$ and
$$
\|\psi f\|_{\sigma_{2},p}=\frac{\sigma_{2}(y)}{\sigma_{1}(y)}|\psi(y)|.
$$
As $y\in T$ was arbitrary, we see that
$$
\left\|\psi/\sigma_{1}\right\|_{\sigma_{2},\infty}=\sup_{y\in T}\frac{\sigma_{2}(y)}{\sigma_{1}(y)}|\psi(y)|\leq \|M_{\psi}\|,
$$ that is, $\psi/\sigma_{1}\in \mathbb{T}_{\sigma_{2},\infty}.$

Conversely, suppose that $\psi/\sigma_{1}\in \mathbb{T}_{\sigma_{2},\infty}.$
For any $f\in \mathbb{T}_{\sigma_{1},p}$ and $n\in \N_{0}$, we easily have
$$
M_{\sigma_{2},p}(n,\psi f)\leq \sup_{y\in T}\left(\frac{\sigma_{2}(y)}{\sigma_{1}(y)}|\psi(y)|\right) M_{\sigma_{1},p}(n,f).
$$
Therefore, $\|\psi f\|_{\sigma_{2},p}\leq \|\psi/\sigma_{1}\|_{\sigma_{2},\infty}\,\|f\|_{\sigma_{1},p}$ for all $f\in \mathbb{T}_{\sigma_{1},p}.$
This gives that, $M_{\psi}:\mathbb{T}_{\sigma_{1},p}\rightarrow \mathbb{T}_{\sigma_{2},p}$
 is a bounded operator. From this proof, it is also clear that
 $$\|M_{\psi}\|=\left\|\psi/\sigma_{1}\right\|_{\sigma_{2},\infty}.
 $$
\epf

By taking $\psi(x)=1$ in the previous theorem, we have the following interesting result.
\bcor\label{incl}
For $0< p\leq\infty$, define the identity operator $I:\mathbb{T}_{\sigma_{1},p}\rightarrow \mathbb{T}_{\sigma_{2},p}$ by $I(f)=f$. Then, $I$
 is a bounded operator $($or equivalently, $\mathbb{T}_{\sigma_{1},p} \subseteq\mathbb{T}_{\sigma_{2},p})$ if and only if
$$\sup_{x\in T}\frac{\sigma_{2}(x)}{\sigma_{1}(x)}<\infty.
$$
\ecor

\section{Compact weighted composition operators}\label{sec:cmpt}
This section is devoted to compactness of weighted composition operators. For any function
 $f$ on $T$, supremum of $|f|$ over an empty set is regarded as zero, that is,
$$\sup_{x\in \emptyset}|f(x)|=0.
$$
We now begin the section with the following characterization of compact weighted
composition operators between $\mathbb{T}_{\sigma, p}$ spaces.

\bthm \label{cmpt}
Let $0<p,q\leq \infty$. Then, the weighted composition operator
$W_{\psi,\varphi}:\mathbb{T}_{\sigma_{1}, p}\rightarrow \mathbb{T}_{\sigma_{2},q}$
 is compact if and only if for every bounded sequence $\{f_{n}\}$ in
 $\mathbb{T}_{\sigma_{1},p}$ that converges to $0$ pointwise, the sequence $\|W_{\psi,\varphi}(f_{n})\|_{\sigma_{2},q}\rightarrow 0$ as $n \rightarrow\infty.$
\ethm
Proof is similar to that of \cite[Theorem 6]{MSha1} and so we omit its proof.

\bthm
The operator $W_{\psi,\varphi}:\mathbb{T}_{\sigma_{1},\infty}\rightarrow \mathbb{T}_{\sigma_{2},\infty}$ is compact if and only if
$$
\sup_{x\in \varphi^{-1}(y)}\left|\frac{\sigma_{2} \psi}{\sigma_{1}\circ \varphi}(x)\right |\rightarrow 0 \; \; \text{ as }  \; \; |y|\rightarrow\infty.
$$
\ethm

\bpf
Suppose that $W_{\psi,\varphi}:\mathbb{T}_{\sigma_{1},\infty}\rightarrow \mathbb{T}_{\sigma_{2},\infty}$ is compact.
For $y\in T$, define $f_{y}$ as
$f_{y}=\frac{1}{\sigma_{1}}\chi_{y}.$ Then,
\begin{align*}
\|f_{y}\|_{\sigma_{1},\infty}= \sup_{x\in T}|\chi_{y}(x)|=1 \text{ for all } y\in T.
\end{align*}
It is obvious that $f_{y}\rightarrow 0$ as $|y|\rightarrow\infty$ and
\begin{align*}
\|W_{\psi,\varphi}(f_{y})\|_{\sigma_{2},\infty} & = \sup_{x\in T}|\sigma_{2}(x) \psi(x) f_{y}(\varphi(x))|
  = \sup_{x\in \varphi^{-1}(y)}\left|\frac{\sigma_{2} \psi}{\sigma_{1}\circ \varphi}(x)\right|.
 \end{align*}
As $W_{\psi,\varphi}$ is compact, by Theorem \ref{cmpt}, we get that
$$\sup_{x\in \varphi^{-1}(y)}\left|\frac{\sigma_{2} \psi}{\sigma_{1} \circ \varphi}(x)\right|\rightarrow 0 \; \; \text{ as }  \; \; |y|\rightarrow\infty.
$$

Conversely, assume that \begin{align*}
\sup_{x\in \varphi^{-1}(y)}\left|\frac{\sigma_{2} \psi}{\sigma_{1}\circ \varphi}(x)\right|\rightarrow 0 \; \; \text{ as }  \; \; |y|\rightarrow\infty.
\end{align*}
Let $\epsilon > 0$ be given. Then, there exists an $M > 0$ such that
\begin{align*}
\sup_{x\in \varphi^{-1}(y)}\left|\frac{\sigma_{2} \psi}{\sigma_{1}\circ \varphi}(x)\right| < \frac{\epsilon}{2}   \hspace{5pt}\text{whenever} \hspace{5pt} |y|>M.\end{align*}
Equivalently,
\begin{align*}
  \sup_{|\varphi(x)|>M}\left|\frac{\sigma_{2} \psi}{\sigma_{1}\circ\varphi}(x)\right|\leq \frac{\epsilon}{2}.
\end{align*}
Let ${f_{n}}$ be a sequence in $\mathbb{T}_{\sigma_{1},\infty}$ such that
$\|f_{n}\|_{\sigma_{1},\infty}\leq 1$ for all $n\in \mathbb{N}$ and $f_{n}\rightarrow 0$
pointwise as $n \rightarrow\infty.$ As $\{y\in T:|y|\leq M\}$ is a finite set, there
exists an $N\in \mathbb{N}$ such that
\begin{align*}
\sup_{|y|\leq M}|f_{n}(y)\sigma_{1}(y)|< \frac{\epsilon}{2\|W_{\psi,\varphi}\|} \hspace{5pt}
\text{for all} \hspace{5pt}n>N.
\end{align*}
For $n>N$,
\begin{align*}
\sup_{|\varphi(x)|\leq M}|\sigma_{2}(x)\psi(x)f_{n}(\varphi(x))| &
\leq \sup_{|\varphi(x)|\leq M}\left|\frac{\sigma_{2}  \psi}{\sigma_{1}\circ \varphi}(x)\right|\sup_{|\varphi(x)|\leq M}|f_{n}(\varphi(x))\sigma_{1}(\varphi(x))|\\
 &<\|W_{\psi,\varphi}\| \frac{\epsilon}{2\|W_{\psi,\varphi}\|}   = \frac{\epsilon}{2},
\end{align*}
and \begin{align*}
 \sup_{|\varphi(x)|> M}|\sigma_{2}(x)\psi(x)f_{n}(\varphi(x))|
 & \leq \sup_{|\varphi(x)|> M}\left|\frac{\sigma_{2} \psi}{\sigma_{1}\circ \varphi}(x)\right|\sup_{|\varphi(x)|> M}|f_{n}(\varphi(x))\sigma_{1}(\varphi(x))|\\
 & \leq \frac{\epsilon}{2} \|f_{n}\|_{\sigma_{1},\infty}\leq \frac{\epsilon}{2}.
\end{align*}
Therefore,
\begin{align*}
\|W_{\psi,\varphi}(f_{n})\|_{\sigma_2,\infty} = \sup_{x\in T}| \sigma_{2}(x) \psi(x) f_{n}(\varphi(x))| < \epsilon\hspace{5pt} \text{for all} \hspace{5pt} n>N.
\end{align*}
Consequently, $\|W_{\psi,\varphi}(f_{n})\|_{\sigma_2,\infty} \rightarrow 0 $ as $n \rightarrow\infty.$ Hence, $W_{\psi,\varphi}$ is compact by Theorem \ref{cmpt}.
\epf

Though we have characterization for compactness of $C_\varphi$ and $M_\psi$ by above
theorem,
we state it for $M_\psi$ because of its independent interest and simplicity.
\bcor\label{cmpt-cor}
 The multiplication operator
 $M_{\psi}:\mathbb{T}_{\sigma_{1},\infty}\rightarrow \mathbb{T}_{\sigma_{2},\infty}$
 is compact if and only if
$$
\frac{\sigma_{2}(x)}{\sigma_{1}(x)}|\psi(x)|\rightarrow 0 \; \mbox{ as }
 |x|\rightarrow\infty.
 $$
\ecor

\bthm
The operator $W_{\psi,\varphi}:\mathbb{T}_{\sigma_{1},p}\rightarrow \mathbb{T}_{\sigma_{2},\infty}$ is compact if and only if
$$
\sup_{x\in\varphi^{-1}(y)}\frac{\sigma_{2}(x)|\psi(x)|}{\sigma_{1}(\varphi(x))}
\left(c_{|\varphi(x)|}\right)^{1/p}
\rightarrow 0 \; \; \mbox{ as }  \; \; |y|\rightarrow\infty.
$$
\ethm

\bpf
Suppose that $W_{\psi,\varphi} : \mathbb{T}_{\sigma_{1},p}\rightarrow \mathbb{T}_{\sigma_{2},\infty}$ is compact. For each $y\in T,$ consider a function
$f_y$ defined by
$$f_{y}(x)=\frac{\left(c_{|y|}\right)^{{1}/{p}}}{\sigma_{1}(y)}\chi_{y}(x).$$
Then, it is easy to see that $\|f_{y}\|_{\sigma_{1},p}=1$ for all $y\in T$ and
$f_{y} \rightarrow 0$ pointwise as $|y|\rightarrow\infty.$ Now,
\begin{align*}
\|W_{\psi,\varphi}(f_{y})\|_{\sigma_{2},\infty} & = \sup_{x\in T}|\sigma_{2}(x)\psi(x)f_{y}(\varphi(x))| \\
& = \sup_{x\in T}\sigma_{2}(x)|\psi(x)| \frac{\left(c_{|y|}\right)^{{1}/{p}}}{\sigma_{1}(y)}\chi_{y}(\varphi(x))\\
 & = \sup_{x\in \varphi^{-1}(y)}\frac{\sigma_{2}(x)|\psi(x)|}{\sigma_{1}
 (\varphi(x))} c_{|\varphi(x)|}^{ {1}/{p}}.
\end{align*}
Thus, the desired result follows from Theorem \ref{cmpt}.

Conversely, assume that
$$
\sup_{x\in \varphi^{-1}(y)}\frac{\sigma_{2}(x)|\psi(x)|}{\sigma_{1}(\varphi(x))}c_
{|\varphi(x)|}^{ {1}/{p}} \rightarrow 0  \; \; \text{ as }  \; \; |y|\rightarrow\infty.
$$
Then, for a given $\epsilon>0$, there exists an $M>0$ such that
$$
\sup_{|\varphi(x)|>M} \frac{\sigma_{2}(x)|\psi(x)|}{\sigma_{1}(\varphi(x))}c_{|\varphi(x)|}
^{ {1}/{p}} \leq\frac{\epsilon}{2}.
$$
Choose a sequence ${(f_{n})}$  in $\mathbb{T}_{\sigma_{1},p}$ such that $\|f_{n}\|_{\sigma_{1},p}\leq1$ for all $n\in \mathbb{N}$ and $f_{n}\rightarrow0$ pointwise as $n\rightarrow\infty.$ As $\{y\in T:|y|\leq M\}$ is a finite set, there exists an
$N\in \mathbb{N}$ such that
$$
\sup_{|y|\leq M}\frac{|f_{n}(y)|\sigma_{1}(y)}{\left(c_{|y|}\right)^{ {1}/{p}}}< \frac{\epsilon}{2\|W_{\psi,\varphi}\|} \; \; \text{ for } \; \; n>N.
$$
For any $n>N$, we have
\begin{align*}
 \sup_{|\varphi(x)|\leq M}|\sigma_{2}(x)\psi(x)f_{n}(\varphi(x))| & \leq \sup_{|\varphi(x)|\leq M}\frac{\sigma_{2}(x)|\psi(x)|}{\sigma_{1}(\varphi(x))}c_{|\varphi(x)|}^{{1}/{p}}
 \sup_{|\varphi(x)|\leq M}\frac{|f_{n}(\varphi(x))|\sigma_{1}(\varphi(x))}{c_{|\varphi(x)|}
 ^{ {1}/{p}}} \\
& <\|W_{\psi,\varphi}\| \frac{\epsilon}{2\|W_{\psi,\varphi}\|}= \frac{\epsilon}{2},
\end{align*}
 and by growth estimate, we get
  \begin{align*}
 \sup_{|\varphi(x)| > M}|\sigma_{2}(x)\psi(x)f_{n}(\varphi(x))|
 & \leq \sup_{|\varphi(x)| > M}
 \frac{\sigma_{2}(x)|\psi(x)|}{\sigma_{1}(\varphi(x))}c_{|\varphi(x)|}^{{1}/{p}}
 \|f_{n}\|_{\sigma_{1},p}
\leq \frac{\epsilon}{2} \|f_{n}\|_{\sigma_{1},p}\leq \frac{\epsilon}{2}.
\end{align*}
Thus,
$$
\|W_{\psi,\varphi}f_{n}\|_{\sigma_{2},\infty}= \sup_{x\in T}|\sigma_{2}(x)\psi(x)f_{n}(\varphi(x))|<\epsilon
$$ for all $n>N$. This gives that
  $\|W_{\psi,\varphi}f_{n}\|\rightarrow 0$ as $n\rightarrow\infty$. Therefore,
  $W_{\psi,\varphi}$ is a compact operator by Theorem \ref{cmpt}.
\epf

As an immediate consequence, we have the following result.
\bcor
For $0< p<\infty$, the multiplication operator $M_{\psi}:\mathbb{T}_{\sigma_{1},p}\rightarrow \mathbb{T}_{\sigma_{2},\infty}$  is compact if and only if
$$
\frac{\sigma_{2}(x)}{\sigma_{1}(x)}|\psi(x)|\left(c_{|x|}\right)^{1/p}\rightarrow 0
\; \mbox{ as } |x|\rightarrow\infty.
$$
\ecor

\bthm
For $0< p<\infty$, the operator $W_{\psi, \varphi}:\mathbb{T}_{\sigma_{1},\infty}\rightarrow \mathbb{T}_{\sigma_{2},p}$ is compact if and only if
$$\sup_{n\in \mathbb{N}_{0}}\frac{1}{c_{n}}\sum_{\substack
 {|x|=n \\ |\varphi(x)|\geq m }} \left|\frac{\sigma_{2} \psi}{\sigma_{1}\circ\varphi}(x)\right|^{p}\rightarrow 0 \; \text{ as }
 m\rightarrow\infty.$$
\ethm
\bpf
Suppose $W_{\psi, \varphi}$ is compact.
For each $m\in \mathbb{N}$, define $g_{m}$ as follows:
$$ g_{m}(x)=\left\{\begin{array}{cl} 0 & \mbox {if} \; |x|<m, \\
 \frac{1}{\sigma_1(x)} & \mbox {if} \; |x|\geq m. \end{array} \right.
 $$
Then, $\|g_{m}\|_{\sigma_{1},\infty}=1$ for all $m\in \mathbb{N}$ and $g_{m}
\rightarrow 0$ pointwise as $m \rightarrow\infty.$\\
For $n\in \mathbb{N}_{0}$,
$$M_{\sigma_{2},p}^{p}(n,W_{\psi,\varphi}(g_{m}))=\frac{1}{c_{n}}\sum_{\substack
 {|x|=n \\ |\varphi(x)|\geq m }}
\frac{\sigma_{2}(x)^{p}|\psi(x)|^{p}}{\sigma_{1}^{p}(\varphi(x))}.
$$
Hence,$$
\|W_{\psi, \varphi}(g_{m})\|=\sup_{n\in \mathbb{N}_{0}}\frac{1}{c_{n}}
\sum_{\substack
 {|x|=n \\ |\varphi(x)|\geq m }}\left|\frac{\sigma_{2}\psi}{\sigma_{1} \circ \varphi}(x)\right|^{p}\rightarrow 0 \text { as } m\rightarrow\infty,$$
 by Theorem \ref{cmpt}.

Conversely, suppose that
\begin{equation*}
\displaystyle \sup_{n\in \mathbb{N}_{0}}\frac{1}{c_{n}}\sum_{\substack
 {|x|=n \\ |\varphi(x)|\geq m }}\left|\frac{\sigma_{2}\psi}{\sigma_{1} \circ \varphi}(x)\right|^{p}\rightarrow 0 \text{ as } m \rightarrow\infty.
\end{equation*}

Let $\epsilon>0$ be given. Then, there exists an $M\in \mathbb{N}$ such that
\begin{equation*}
\displaystyle \sup_{n\in \mathbb{N}_{0}}\frac{1}{c_{n}}\sum_{\substack
 {|x|=n \\ |\varphi(x)|\geq M }}\left|\frac{\sigma_{2}\psi}{\sigma_{1} \circ \varphi}(x)\right|^{p}<\frac{\epsilon}{2}.
 \end{equation*}
Let $\{f_{k}\}$ be  a sequence in the unit ball of $\mathbb{T}_{\sigma_{1},\infty}$ such that $f_{k}\rightarrow 0$ pointwise as $k \rightarrow\infty.$ Since $\{y\in T: |y|<M\}$ is a finite set, there exists an $N\in \mathbb{N}$ such that
\begin{equation*}
\displaystyle \sup_{|y|<M}|f_{k}(y)|^{p}\sigma_{1}^{p}(y)<\frac{\epsilon}
{2\|W_{\psi,\varphi}\|^{p}} \text{ for all } k>N.
\end{equation*}
Thus, for $n\in \mathbb{N}_{0}$ and $k>N$,
\beqq
M_{\sigma_{2},p}^{p}(n, W_{\psi,\varphi}f_{k})&=&\frac{1}{c_{n}}\sum_{\substack
 {|x|=n \\ |\varphi(x)|\geq M}}\sigma_{2}^{p}(x)|\psi(x)|^{p}|f_{k}(\varphi(x))|^{p}+\\
& &\frac{1}{c_{n}}\sum_{\substack
 {|x|=n \\ |\varphi(x)|<M}}\sigma_{2}^{p}(x)|\psi(x)|^{p}|f_{k}(\varphi(x))|^{p}.
\eeqq
Therefore, for all $k>N$,
\beqq
\|W_{\psi,\varphi}f_{k}\|_{\sigma_{2},p}^{p}
&\leq& \sup_{n\in \mathbb{N}_{0}}\frac{1}{c_{n}}\sum_{\substack
 {|x|=n \\ |\varphi(x)|\geq M }}\left|\frac{\sigma_{2}\psi}{\sigma_{1} \circ \varphi}(x)\right|^{p}|f_{k}(\varphi(x))|^{p}\sigma_{1}^{p}(\varphi(x))+\\
& & \sup_{n\in \mathbb{N}_{0}}\frac{1}{c_{n}}\sum_{\substack
 {|x|=n \\ |\varphi(x)|<M}}\left|\frac{\sigma_{2}\psi}{\sigma_{1} \circ \varphi}(x)\right|^{p}|f_{k}(\varphi(x))|^{p}\sigma_{1}^{p}(\varphi(x))\\
&\leq& \|f_{k}\|_{\sigma_{1},\infty}^{p} \sup_{n\in \mathbb{N}_{0}}\frac{1}{c_{n}}
\sum_{\substack{|x|=n \\ |\varphi(x)|\geq M }}\left|\frac{\sigma_{2}\psi}{\sigma_{1}
\circ \varphi}(x)\right|^{p}+
 \|W_{\psi,\varphi}\|^{p} \sup_{|\varphi(x)|<M}|f_{k}(\varphi(x))|^{p}\sigma_{1}^{p}(\varphi(x))\\
&<& \frac{\epsilon}{2}+\|W_{\psi,\varphi}\|^{p}\frac{\epsilon}{2\|W_{\psi,\varphi}\|^{p}}
 =\epsilon.
\eeqq
 Thus, $\|W_{\psi,\varphi}f_{k}\|\rightarrow 0$ as $k \rightarrow\infty $.
Hence, by Theorem \ref{cmpt}, $W_{\psi,\varphi}$ is a compact operator. It completes the proof.
\epf

\bcor
 The multiplication operator $M_{\psi}:\mathbb{T}_{\sigma_{1},\infty}\rightarrow \mathbb{T}_{\sigma_{2},p}$, $p<\infty$ is compact if and only if
$$
M_{\sigma_{2},p}\left(n,\frac{\psi}{\sigma_{1}}\right)\rightarrow 0  \; \mbox{ as } |x|\rightarrow\infty.
$$
\ecor

\bthm
Consider the operator $W_{\psi,\varphi}:\mathbb{T}_{\sigma_{1},p}\rightarrow \mathbb{T}_{\sigma_{2},p}$, $0<p<\infty$. Then,
\begin{itemize}
  \item[(i)] the operator $W_{\psi,\varphi}$ is compact if
\begin{equation*}
\displaystyle \sup_{n\in \mathbb{N}_{0}}\frac{1}{c_{n}}\sum_{\substack
 {|x|=n \\ |\varphi(x)|\geq m }}\left|\frac{\sigma_{2}\psi}{\sigma_{1} \circ \varphi}(x)\right|^{p}c_{|\varphi(x)|}\rightarrow 0 \text{ as } m \rightarrow\infty.
\end{equation*}

  \item[(ii)]if $W_{\psi,\varphi}$ is a compact operator, then
  $$\sup_{n\in \mathbb{N}_{0}}\frac{1}{c_{n}}\sum_{\substack
 {|x|=n \\ |\varphi(x)|\geq m }}\left|\frac{\sigma_{2}\psi}{\sigma_{1} \circ \varphi}(x)\right|^{p}\rightarrow 0 \text{ as } m \rightarrow\infty.
 $$
\end{itemize}
\ethm

\bpf Consider the operator $W_{\psi,\varphi}:\mathbb{T}_{\sigma_{1},p}\rightarrow \mathbb{T}_{\sigma_{2},p}$, $0<p<\infty$.
Note that the  function $1/\sigma_1\in \mathbb{T}_{\sigma_{1},p}$ with $\|1/\sigma_1\|_{\sigma_{1},p}=1.$ This gives that $\psi/(\sigma_{1}\circ\varphi)\in \mathbb{T}_{\sigma_{2},p}$ and $\|\psi/(\sigma_{1}\circ\varphi)\|_{\sigma_{2},p}\leq \|W_{\psi,\varphi}\|$, i.e.,
\begin{equation*}
\displaystyle \sup_{n\in \mathbb{N}_{0}}\frac{1}{c_{n}}\sum_{|x|=n}\left|\frac{\sigma_{2}\psi}{\sigma_{1} \circ \varphi}(x)\right|^{p}\leq \|W_{\psi,\varphi}\|^{p}.
\end{equation*}

Now, we prove $\textrm{(i)}$. Assume that
\begin{equation*}
\displaystyle \sup_{n\in \mathbb{N}_{0}}\frac{1}{c_{n}}\sum_{\substack
 {|x|=n \\ |\varphi(x)|\geq m }}\left|\frac{\sigma_{2}\psi}{\sigma_{1} \circ \varphi}(x)\right|^{p}c_{|\varphi(x)|}\rightarrow 0 \mbox{ as } m \rightarrow\infty.
\end{equation*}
For a given $\epsilon>0,$ there exists an $M\in \mathbb{N}$ such that
\begin{equation*}
\displaystyle \sup_{n\in \mathbb{N}_{0}}\frac{1}{c_{n}}\sum_{\substack
 {|x|=n \\ |\varphi(x)|\geq M }}\left|\frac{\sigma_{2}\psi}{\sigma_{1} \circ \varphi}(x)\right|^{p}c_{|\varphi(x)|}<\frac{\epsilon}{2}.
\end{equation*}
Let ${f_{k}}$ be a sequence in $\mathbb{T}_{\sigma_{2},p}$ such that $\|f_{k}\|_
{\sigma_{1},p}\leq1$ for all $k\in \mathbb{N}$ and $f_{k}\rightarrow 0$ pointwise as $k \rightarrow\infty.$
As $\{y\in T:|y|<M\}$ is a finite set, there exists an $N\in \mathbb{N}$ such that
\begin{equation*}
\displaystyle \sup_{|y|<M}|f_{k}(y)|^{p}\sigma_{1}^{p}(y)<\frac{\epsilon}{2\|W_{\psi,\varphi}\|^{p}}
\end{equation*}
for all $k>N$.
Thus, for any $k>N$ and $n\in \mathbb{N}_{0}$, we have
\beqq
M_{\sigma_{2},p}^{p}(n, W_{\psi,\varphi}f_{k}) &=&\frac{1}{c_{n}}\sum_{\substack
 {|x|=n \\ |\varphi(x)|\geq M }}\sigma_{2}^{p}(x)|\psi(x)|^{p}|f_{k}(\varphi(x))|^{p}+\\
 & &
\frac{1}{c_{n}}\sum_{\substack
 {|x|=n \\ |\varphi(x)|< M}}\sigma_{2}^{p}(x)|\psi(x)|^{p}|f_{k}(\varphi(x))|^{p}.
\eeqq

Therefore, for all $k>N$,
\beqq
\|W_{\psi,\varphi}f_{k}\|_{\sigma_{2},p}^{p} &\leq& \sup_{n\in \mathbb{N}_{0}}\frac{1}{c_{n}}\sum_{\substack
 {|x|=n \\ |\varphi(x)|<M }}\left|\frac{\sigma_{2}\psi}{\sigma_{1} \circ \varphi}(x) \right|^{p}|f_{k}(\varphi(x))|^{p}\sigma_{1}^{p}(\varphi(x))+ \\
 & & \sup_{n\in \mathbb{N}_{0}}\frac{1}{c_{n}}\sum_{\substack
 {|x|=n \\ |\varphi(x)|\geq M }}\left|\frac{\sigma_{2}\psi}{\sigma_{1} \circ \varphi}(x) \right|^{p}c_{|\varphi(x)|}\frac{|f_{k}(\varphi(x))|^{p}\sigma_{1}^{p}(\varphi(x))}{c_{|\varphi(x)|}} \\
 &\leq& \|W_{\psi,\varphi}\|^{p} \sup_{|\varphi(x)|<M}|f_{k}(\varphi(x))|^{p} \sigma_{1}^{p}(\varphi(x)) +\\ & &
\|f_{k}\|^{p}_{\sigma_{2},p} \sup_{n\in \mathbb{N}_{0}}\frac{1}{c_{n}}\sum_{\substack
 {|x|=n \\ |\varphi(x)|\geq M }}\left|\frac{\sigma_{2}\psi}{\sigma_{1} \circ \varphi}(x) \right|^{p}c_{|\varphi(x)|}
\\
&<&\|W_{\psi,\varphi}\|^{p} \frac{\epsilon}{2\|W_{\psi,\varphi}\|^{p}}+ \frac{\epsilon}{2} =\epsilon.
\eeqq
This implies that $\|W_{\psi,\varphi}f_{k}\|\rightarrow 0 \text{ as } k \rightarrow\infty.$
Hence, by Theorem \ref{cmpt}, $W_{\psi, \varphi}$ is a compact operator.

Next, we prove $\textrm{(ii)}$. Suppose that $W_{\psi, \varphi}$ is a compact operator.
For each $m\in \mathbb{N}$, define
$$
g_{m}(x)=\left\{\begin{array}{cl} 0 & \mbox {if} \; |x|<m, \\
\frac{1}{\sigma_1(x)} & \mbox {if} \; |x|\geq m. \end{array} \right.
$$
Then, $\|g_{m}\|_{\sigma_{1},p}=1$ for all $m\in \mathbb{N}$ and $g_{m}
\rightarrow 0$ pointwise as $m \rightarrow\infty$.
Hence,
$$\|W_{\psi, \varphi}(g_{m})\|=\sup_{n\in \mathbb{N}_{0}}\frac{1}{c_{n}}
\sum_{\substack
 {|x|=n \\ |\varphi(x)|\geq m }}\left|\frac{\sigma_{2}\psi}{\sigma_{1} \circ \varphi}(x)\right|^{p}\rightarrow 0$$
  as $m \rightarrow\infty$, by Theorem \ref{cmpt}. It completes the proof.
\epf

\bthm\label{cmpt-MO}
For $0<p\leq\infty$, the multiplication operator $M_{\psi}:\mathbb{T}_{\sigma_{1},p}\rightarrow \mathbb{T}_{\sigma_{2},p}$  is compact if and only if $$\frac{\sigma_{2}(x)}{\sigma_{1}(x)}|\psi(x)|\rightarrow 0 \; \mbox{  as } |x|\rightarrow\infty.$$
\ethm
\bpf
For the case $p=\infty$, refer Corollary \ref{cmpt-cor}. So, it is enough to consider
the case $p<\infty$.

Suppose $M_{\psi}$ is compact.
For each $y\in T$, take $f_y$ as
$$ f_{y}(x) = \left\{
 \begin{array}{cl}
   \frac{\left(c_{|y|}\right)^{{1}/{p}}}{\sigma_{1}(y)}  & \mbox{ if }  x=y,\\
0  & \mbox{ otherwise}.
\end{array}
\right.
$$
It is evident that $\|f_y\|_{\sigma_{1},p}=1$ for every $y\in T$.
As $f_{y}\rightarrow 0$ pointwise as $|y|\rightarrow\infty$, we have
$$\|M_{\psi}(f_{y})\|=\frac{\sigma_{2}(y)}{\sigma_{1}(y)}|\psi(y)|\rightarrow 0
 \; \mbox{  as } |y|\rightarrow\infty.
 $$

Conversely, suppose that
$$\frac{\sigma_{2}(x)}{\sigma_{1}(x)}|\psi(x)|\rightarrow 0 \mbox{ as } |x|\rightarrow\infty.$$
For each $n\in \mathbb{N}$, consider the map $$ \psi_{n}(x) = \left\{
 \begin{array}{cl}
  \psi(x) & \mbox{ if }  |x|<n,\\
0  & \mbox{ otherwise}.
\end{array}
\right.
$$ 
Thus, $M_{\psi_{n}}$ is a finite  rank operator and hence compact. Further,
\begin{align*}
\|M_{\psi}-M_{\psi_{n}}\| & =\|M_{\psi-\psi_{n}}\|=\sup_{x\in T}\frac{\sigma_{2}(x)}{\sigma_{1}(x)}|(\psi-\psi_{n})(x)|\\ & = \sup_{|x|\geq n}\frac{\sigma_{2}(x)}{\sigma_{1}(x)}
|\psi(x)| \rightarrow 0\; \mbox{ as } n \rightarrow\infty.
 \end{align*}
As $M_{\psi}$ is a limit of a sequence of compact operators, $M_{\psi}$ is also a
compact operator.
\epf

\bcor
For $0<p\leq\infty$, the identity operator $I:\mathbb{T}_{\sigma_{1},p}\rightarrow \mathbb{T}_
{\sigma_{2},p}$ is compact  if and only if
$\frac{\sigma_{2}(x)}{\sigma_{1}(x)}\rightarrow 0 \hspace{5pt} as \hspace{5pt}|x|\rightarrow\infty$.
\ecor

\section{Isometric weighted composition operators}\label{sec:iso}
In this section, we will discuss about isometric weighted
composition operators between $\mathbb{T}_{\sigma, p}$ spaces. Isometric multiplication
and isometric composition operators on $\mathbb{T}_{p}$ spaces were characterized in
\cite[Theorem 4.10]{MP-Tp-spaces} and \cite[Theorems 13 and 14]{CO-Tp-spaces2}, respectively.
When $p\neq q$, there are no isometric multiplication operators from $\mathbb{T}_{p}$ to $\mathbb{T}_{q}$ (see \cite[Theorem 5]{MSha1}).
\bthm
The weighted composition operator $W_{\psi,\varphi}:\mathbb{T}_{\sigma_{1},\infty}
\rightarrow\mathbb{T}_{\sigma_{2},\infty}$ is isometry if and only if $\varphi$ is onto and
$$
\sup_{x\in \varphi^{-1}(y)} \left|\frac{\sigma_{2} \psi}{\sigma_{1}\circ\varphi}(x) \right| = 1\; \text{ for all } \; y\in T.
$$
\ethm
\bpf
Suppose that $W_{\psi,\varphi}:\mathbb{T}_{\sigma_{1},\infty}\rightarrow \mathbb{T}_{\sigma_{2},\infty}$ is an isometry.
We first show that $\varphi$ is onto.
For the contrary, assume that $\varphi$ is not onto. Choose $y\in T\setminus\varphi(T)$.
Then, $f=\chi_{y}\neq0,$ but $\psi(f\circ\varphi)\equiv0$. Thus, $W_{\psi,\varphi}:\mathbb{T}_{\sigma_{1},\infty}\rightarrow \mathbb{T}_{\sigma_{2},\infty}$ cannot be isometry. This contradiction yields that $\varphi$ is onto.

Next, for $y\in T$, take $f_{y}=\frac{1}{\sigma_{1}}\chi_{y}$. Then, $\|f_{y}\|_{\sigma_{1},\infty}=1$ for all $y\in T$.
\begin{align*}
\|\psi(f_{y} \circ \varphi)\|_{\sigma_{2},\infty} & = \sup_{x\in T}\frac{\sigma_{2}(x)|\psi(x)|}{\sigma_{1}(\varphi(x))}\chi_{y}(\psi(x))
=\sup_{x\in \varphi^{-1}(y)}\left|\frac{\sigma_{2}\psi}{\sigma_{1}\circ\varphi}(x)\right|.
\end{align*}
As $W_{\psi,\varphi}$ is an isometry, we see that
$$\sup_{x\in \varphi^{-1}(y)}\left|\frac{\sigma_{2} \psi}{\sigma_{1} \circ \varphi}(x)\right|=1 \; \text{ for all } y\in T.$$
Conversely, assume that $\varphi$ is onto and $$\sup_{x\in \varphi^{-1}(y)}\left|\frac{\sigma_{2} \psi}{\sigma_{1} \circ \varphi}(x)\right|=1 \; \text{ for all } y\in T.$$
As $\varphi$ is onto, $T$ can be written as disjoint union of the sets
 $\varphi^{-1}(y),\, y\in T.$ For any $f\in \mathbb{T}_{\sigma_{1},\infty}$,
\begin{align*}
\|\psi(f\circ\varphi)\|_{\sigma_{2},\infty} & = \sup_{x\in T}\frac{\sigma_{2}(x)|\psi(x)|}{\sigma_{1}(\varphi(x))}|f(\varphi(x))|\sigma_{1}(\varphi(x))\\
& = \sup_{y\in T}\sup_{x\in \varphi^{-1}(y)}\frac{\sigma_{2}(x)|\psi(x)|}{\sigma_{1}(\varphi(x))}|f(y)|\sigma_{1}(y)\\
& = \sup_{y\in T}|f(y)|\sigma(y)
 =\|f\|_{\sigma_{1},\infty}.
\end{align*}
Therefore, $W_{\psi,\varphi}$ is an isometry. Hence the theorem.
 \epf
\bthm
For $0<p<\infty$, there are no isometric weighted composition operators $W_{\psi,\varphi}:\mathbb{T}_{\sigma_{1},p}\rightarrow \mathbb{T}_{\sigma_{2},\infty}$.
\ethm
\bpf
Suppose there exists an isometric operator  $W_{\psi,\varphi}:\mathbb{T}_{\sigma_{1},p}
\rightarrow \mathbb{T}_{\sigma_{2},\infty}$. By using same arguments as in previous theorem,
we get $\varphi$ is onto.
For  $y\in T$, consider
$$f_{y}(x)=\frac{\left(c_{|y|}\right)^{1/p}}{\sigma_{1}(x)}\chi_{y}(x).$$
Then, $\|f_{y}\|_{\sigma_{1},p}=1$ for all $y\in T$. As $W_{\psi,\varphi}$ is an isometry,
 we have
$$\|\psi(f_{y}\circ\varphi)\|_{\sigma_{2},\infty}= \sup_{x\in \varphi^{-1}(y)}\frac{\sigma_{2}(x)|\psi(x)|}{\sigma_{1}
(\varphi(x))}c_{|\varphi(x)|}^{1/p}=1$$ for all $y\in T$.
Now for any $f\in \mathbb{T}_{\sigma,p},$
\begin{align*}
\|\psi(f\circ\varphi)\|_{\sigma_{2},\infty} & = \sup_{x\in T}\frac{\sigma_{2}(x)|\psi(x)|}{\sigma_{1}(\varphi(x))}c_{|\varphi(x)|}^{1/p}
\frac{|f(\varphi(x)|\sigma_{1}(\varphi(x))}{c_{|\varphi(x)|}^{1/p}}\\
& = \sup_{y\in T}\left(\sup_{x\in \varphi^{-1}(y)}
\frac{\sigma_{2}(x)|\psi(x)|}{\sigma_{1}(\varphi(x))}c_{|\varphi(x)|}^{1/p}\right)
\frac{|f(y)|\sigma_{1}(y)}{\left(c_{|y|}\right)^{1/p}}\\
& = \sup_{y\in T}\frac{|f(y)|\sigma_{1}(y)}{\left(c_{|y|}\right)^{1/p}}.
\end{align*}

Choose $y_{1}, y_{2}\in T$ with $|y_{1}| = |y_{2}|$ (say $m$). Take $g$ as
\[ g(x) = \left\{
 \begin{array}{cl}
   \frac{1}{2 \sigma_1(x)}   & \mbox{ if }  x = y_1,\\
   \frac{1}{3 \sigma_1(x)}   & \mbox{ if }  x = y_2,\\
0  & \mbox{ otherwise}.
\end{array}
\right.
\]
Then, $$
\|g\|_{\sigma_{1},p}^{p}=M_{\sigma_{1},p}^{p}(m,g)=\frac{1}{c_{m}}
\left(\frac{1}{2^{p}} +\frac{1}{3^{p}}\right)
$$ and
\begin{align*} \|W_{\psi,\varphi} g\|_{\sigma_{2},\infty}^{p}& = \sup_{y\in T}\frac{|g(y)|^{p}\sigma_{1}^{p}(y)}{c_{|y|}}
 = \max\left\{\frac{1}{2^{p}c_{m}},\frac{1}{3^{p}c_{m}}\right\}
 =\frac{1}{2^{p}c_{m}}. \end{align*}
Therefore, $\|g\|_{\sigma_{1},p}\neq\|W_{\psi,\varphi} g\|_{\sigma_{2},\infty}.$
Hence, $W_{\psi,\varphi} $ cannot be an isometry.
\epf

As there are no isometric multiplication operators from $\mathbb{T}_{\infty}$ to $\mathbb{T}_{p}$ when $p<\infty$ (see \cite[Theorem 5]{MSha1}),
we suspect that there are no isometric weighted composition operators $W_{\psi,\varphi}:\mathbb{T}_{\sigma_{1},\infty}\rightarrow \mathbb{T}_{\sigma_{2},p}$
when $p<\infty$.
But, it is  still an open problem.

\section{Examples}\label{sec:eg}
Consider the operators
$M_{\psi}, C_{\varphi}, W_{\psi,\varphi}:\mathbb{T}_{\sigma,p}\rightarrow \mathbb{T}_{\sigma,p}$. From the definition, it is clear that
$$
W_{\psi,\varphi}=M_{\psi}C_{\varphi}, \mbox{~~i.e.,~~} W_{\psi,\varphi}(f)=M_{\psi}(C_{\varphi}(f)),\, f\in\mathbb{T}_{\sigma,p}.
$$
Suppose $M_{\psi}, C_{\varphi}$ are bounded operators on $\mathbb{T}_{\sigma,p}$.
Then, it is well-known that the operator $W_{\psi,\varphi}(=M_{\psi}C_{\varphi})$ is
 also bounded with $\|W_{\psi,\varphi}\|\leq \|M_{\psi}\|\, \|C_{\varphi}\|$.
  Moreover, if at least one of these operators
$M_{\psi}, C_{\varphi}$ is compact, then the operator $W_{\psi,\varphi}$ is
also compact.

 In this section,
we will compare the boundedness and compactness of operators
$M_{\psi}, C_{\varphi}$ with the boundedness and compactness of the operator $W_{\psi,\varphi}$. Throughout the section, we will assume that $0<p<\infty.$

\beg
For each $n\in \mathbb{N}_{0}$, choose $x_{n}\in T$ with $|x_{n}|=n$. Define
\[ \psi(x) = \left\{
 \begin{array}{ll}
   1   & \mbox{ if }  x = x_{n}   \mbox{ for some }  n\in \mathbb N_{0},\\
0  & \mbox{ otherwise}
\end{array}
\right.
\]
and $\varphi(x) = x_{n}$ when $|x| = n$.

Since $\psi$ is a bounded function, $M_{\psi}$ is a bounded operator on
$\mathbb{T}_{p}$.  Consider the function
\[ f(x) = \left\{
 \begin{array}{cl}
  c_{n}^{1/p}   & \mbox{ if }  x = x_{n},  \;  n\in \mathbb \N_{0},\\
0  & \mbox{ otherwise}.
\end{array}
\right.
\]
Then, $\|f\|_p=1$ as $M_{p}(n,f)=1$ for all $n\in \N_{0}$.
For each $n\in \N_{0},$
\begin{align*} M_{p}^{p}(n,C_{\varphi} f) & = \frac{1}{c_{n}}\sum_{|x|=n}|f(\varphi(x))|^{p}  =\frac{1}{c_{n}}\sum_{|x|=n}|f(x_{n})|^{p} = c_{n}.
\end{align*}
Therefore,
$C_{\varphi}f\notin \mathbb{T}_{p}$ as the sequence $(c_{n})$ is unbounded.
Thus, $C_{\varphi}$ is not a bounded operator on $\mathbb{T}_{p}$.

Next, for $f\in \mathbb{T}_{p}$, $n\in \N_{0}$,
\begin{align*} M_{p}^{p}(n, \psi (f \circ\varphi))=  \frac{1}{c_{n}}\sum_{|x|=n}|\psi(x)f(\varphi(x))|^{p}=
 \frac{1}{c_{n}}|\psi(x_{n}) f(\varphi(x_{n}))|^{p}
=\frac{1}{c_{n}}|f(x_{n})|^{p}\leq \|f\|^{p}.
\end{align*}
Thus, $\|\psi(f\circ\varphi)\|_{p} \leq \|f\|_{p}$ for all $f\in \mathbb{T}_{p}.$
Therefore, $W_{\psi, \varphi}$ is a bounded operator on $\mathbb{T}_{p}.$
\eeg

This example shows that $W_{\psi,\varphi}$ can be a bounded operator, though $C_{\varphi}$ is not.

\beg
For each $n\in \mathbb{N}_0$, fix $x_{n}\in T$ with $|x_{n}|=n$.
Define
\[ \psi(x) = \left\{
 \begin{array}{cl}
  c_{n}^{{1}/{p}}   & \mbox{ if }  x=x_{n}, n\in \N_{0},\\
0  & \mbox{ otherwise}
\end{array}
\right.
\]
and $\varphi(x)\equiv\textsl{o}$, where $\textsl{o}$ is the  root of T.
As $\{c_{n}\}_{n\in \N_{0}}$ is unbounded, $M_{\psi}$ is not a bounded operator
on $\mathbb{T}_{p}$. Note that,
 $\varphi\equiv\textsl{o}$ implies $C_{\varphi}=$ evaluation at $\textsl{o}$, which is
 clearly a bounded linear functional.

For $f\in \mathbb{T}_{p}$ and $n\in \N_{0},$
\begin{align*}
M_{p}^{p}(n,\psi(f\circ\varphi)) & =\frac{1}{c_{n}}\sum_{|x|=n}
|\psi(x)|^{p}|f(o)|^{p} =|f(o)|^{p}\leq \|f\|_{p}^{p}.
\end{align*}
Therefore, $\|\psi(f\circ \varphi)\|_{p}\leq \|f\|_{p}$ for all $f\in \mathbb{T}_{p}$,
which in turn implies that $W_{\psi, \varphi}$ is a bounded operator on $\mathbb{T}_{p}$.
\eeg

This example shows that  $W_{\psi,\varphi}$ can be a bounded operator though
$M_{\psi}$ is not a bounded operator.

\beg
As before, for each $n\in\N_{0}$, fix $x_{n}\in T$ with $|x_{n}|=n$.
Choose a subsequence $\{c_{n_{k}}\}$ of $\{c_{n}\}$ such that $c_{n_{k}}
\rightarrow\infty$ as $k\rightarrow\infty$. In particular, $c_{n_{2k+1}}
 \mbox{~and~} c_{n_{2k}}\rightarrow\infty$ as $k\rightarrow\infty$.
Consider the maps
\[ \psi(x) = \left\{
 \begin{array}{cl}
  c_{n}^{{1}/{p}}   & \mbox{ if }  x = x_{n}, n=n_{2k+1},\; k\in \mathbb{N},\\
0  & \mbox{ otherwise}
\end{array}
\right.
\]
and \[ \varphi(x) = \left\{
 \begin{array}{cl}
  x_{n_{2k}}   & \mbox{ if }  |x| = {n_{2k}}, \;  k\in \mathbb N_{0},\\
\textsl{o } & \mbox{ otherwise}.
\end{array}
\right.
\]

Since  $c_{n_{2k+1}} \rightarrow\infty$  as $k \rightarrow\infty,$ we get that
 $\psi$ is an unbounded function and hence $M_{\psi}$ is not a bounded operator
  on $\mathbb{T}_{p}$.
Take \[ f(x) = \left\{
 \begin{array}{cl}
  c_{n}^{{1}/{p}}   & \mbox{ if }  x=x_{n}, \;n\in \N_{0},\\
0  & \mbox{ otherwise}.
\end{array}
\right.
\] 
Then, $\|f\|_{p}=1$.
For each $k\in \mathbb{N}$,
\begin{align*}
M_{p}^{p}(n_{2k},C_{\varphi}f)
& =\frac{1}{c_{n_{2k}}}\sum_{|x|=n_{2k}}|f(\varphi(x))|^{p}\\
& =\frac{1}{c_{n_{2k}}}\sum_{|x|=n_{2k}}|f(x_{n_{2k}})|^{p} =
 c_{n_{2k}} \rightarrow\infty \mbox{ as } k \rightarrow\infty.
\end{align*}
Thus $C_{\varphi}f\notin \mathbb{T}_{p}$, and hence $C_{\varphi}$ is not a bounded
operator on $\mathbb{T}_{p}$.

For $n=n_{2k}$, $k\in \mathbb{N}_0$ and $f\in \mathbb{T}_{p},$ it is obvious that
$M_{p}^{p}(n,\psi(f\circ\varphi))=0,$
as $\psi(x)=0$ for all $x$ with $|x|=n_{2k}.$
For $n\neq n_{2k}, k \in \mathbb{N}_0$,
\begin{align*}
M_{p}^{p}(n, \psi(f\circ \varphi)) &
=\frac{1}{c_{n}}\sum_{|x|=n}|\psi(x)|^{p}|f(o)|^{p}
 = \frac{1}{c_{n}}|\psi(x_{n})|^{p}|f(o)|^{p}=|f(o)|^{p}\leq \|f\|_{p}^{p}.
\end{align*}
Therefore, $\|\psi(f\circ \varphi)\|_{p}\leq \|f\|_{p}$ for all $f\in \mathbb{T}_{p}$.
Hence $W_{\psi, \varphi}$ is a bounded operator on $\mathbb{T}_{p}$.
\eeg

This example shows that $W_{\psi, \varphi}$ can be a bounded operator even though both $M_{\psi}$ and $C_{\varphi}$ are not bounded operators.

 \beg
 Choose the weights $\sigma_{1}, \sigma_{2}$ such that $\sigma_{2}/\sigma_{1}$
 is unbounded.  For example, $\sigma_{1}(x)=|x|+1$ and $\sigma_{2}(x)=(|x|+1)^{2}$.
Take $\varphi(x)=x$ , $x\in T.$
As $\sigma_{2}/\sigma_{1}$ is unbounded, the operator $C_{\varphi}$ is not bounded,
 by Corollary \ref{incl}.
Define $$
\psi(x)=\frac{1}{1+|x|}\left(\frac{\sigma_{1}(x)}{\sigma_{2}(x)}\right),
$$ so that $\frac{\sigma_{2}(x)}{\sigma_{1}(x)}|\psi(x)|\rightarrow 0$ as $|x|\rightarrow\infty.$
Thus, $W_{\psi,\varphi}=M_{\psi}$ is a compact operator from $\mathbb{T}_{\sigma_{1},p}$ to $\mathbb{T}_{\sigma_{2},p}$, by Theorem \ref{cmpt-MO}.
\eeg

The above example shows that the operator $W_{\psi,\varphi}:\mathbb{T}_{\sigma_{1},p}\rightarrow \mathbb{T}_{\sigma_{2},p}$ can be compact, even if $C_{\varphi}$ is an unbounded operator.

\beg
Choose $\sigma_{1}$ such that $1/\sigma_{1}$  is unbounded $($eg. $\sigma_{1}(x)=1/(|x|+1))$. Take $\varphi$ to be a constant map. That is, $\varphi(x)=y$ for all $x\in T$ and for some
$y\in T$.  Define $\psi=1/\sigma_{2}$, so that $\psi\in \mathbb{T}_{\sigma_{2},p}$ with $\|\psi\|_{\sigma_{2},p}=1.$
Consider the operator $M_{\psi}:\mathbb{T}_{\sigma_{1},p}\rightarrow \mathbb{T}_{\sigma_{2},p}.$
As $\sigma_{2}\psi/\sigma_{1}$ is unbounded, $M_{\psi}$  is not a bounded operator by
Theorem \ref{bdd-Mo}.
For any $f\in \mathbb{T}_{\sigma_{1},p}$ ,
\begin{align*}\|W_{\psi,\varphi}f\|_{\sigma_{2},p} & =\|\psi f(y)\|_{\sigma_{2},p} =|f(y)|\,\|\psi\|_{\sigma_{2},p} =|f(y)|\leq
\frac{\left(c_{|y|}\right)^{{1}/{p}}}{\sigma_{1}(y)}\|f\|_{\sigma_{1},p}.
\end{align*}
Thus, $W_{\psi,\varphi}:\mathbb{T}_{\sigma_{1},p}\rightarrow \mathbb{T}_{\sigma_{2},p}$
 is a bounded operator.

 To prove the compactness, start with a bounded sequence $\{f_{n}\}$ in $\mathbb{T}_{\sigma_{1},p}$ such that $f_{n} \rightarrow 0$ pointwise. Then,
$$\|W_{\psi,\varphi}f_{n}\|_{\sigma_{2},p}=|
f_n(y)|\,\|\psi\|_{\sigma_{2},p}=|f_{n}(y)|\rightarrow 0 \; \mbox{ as } n\rightarrow\infty.$$
Therefore, $W_{\psi,\varphi}$ is a compact operator, by Theorem \ref{cmpt}.
\eeg

This example shows that the operator $W_{\psi,\varphi}:\mathbb{T}_{\sigma_{1},p}\rightarrow \mathbb{T}_{\sigma_{2},p}$ can be compact even if $M_{\psi}$ is an unbounded operator.

\beg
For each $n\in\N_{0}$, fix $x_{n}\in T$ with $|x_{n}|=n$. Choose $\sigma_{1}$ such that $\sigma_{1}(x_{n})\rightarrow 0$ as $n \rightarrow\infty\, ($eg. $\sigma_{1}(x_{n})=\frac{1}{n+1})$, and $\sigma_{2}$ such that
$\sum_{|x|=n}\sigma_{2}^{p}(x)\geq1$ if $n$ is even $($eg.
$\sigma_{2}(x)=1$ for all $x\in T)$.
Take $$ \psi(x) = \left\{
 \begin{array}{cl}
  1/\sigma_{2}(x) & \mbox{ if }  x= x_{2k-1}, k\in \mathbb{N},\\
0  & \mbox{ otherwise}.
\end{array}
\right.
$$
Since $$\sup_{x\in T}\frac{\sigma_{2}(x)}{\sigma_{1}(x)}|\psi(x)|= \sup_{k\in \mathbb{N}}\frac{1}{\sigma_{1}(x_{2k-1})}$$ is not finite, $M_{\psi}:\mathbb{T}_{\sigma_{1},p}\rightarrow \mathbb{T}_{\sigma_{2},p}$ is not a bounded operator, by Theorem \ref{bdd-Mo}.

Let us now define,
 $$ \varphi(x) = \left\{
 \begin{array}{cl}
  x_{n} & \mbox{ if }  |x| \mbox{ is even},\\
 \textsl{o}& \mbox{ if }  |x| \mbox{ is odd}.
\end{array}
\right.
$$
Consider the map
 $$ f(x) = \left\{
 \begin{array}{cl}
  \frac{\left(c_{n}\right)^{{1}/{p}}}{\sigma_{1}(x_{n})} & \mbox{ if } x=x_{n}, n\in \N_{0},\\
0  & \mbox{ otherwise}.
\end{array}
\right.
$$
Then, $f\in \mathbb{T}_{\sigma_{1},p}$ with $\|f\|_{\sigma_{1},p}=1$. For $n=2k$, $k\in \mathbb{N}_0$,
\begin{align*}
M_{\sigma_{2},p}(n, C_{\varphi}(f)) & =\frac{1}{c_{n}^{{1}/{p}}}\left(\sum_{|x|=n}\sigma_{2}^{p}(x)\right)^{{1}/{p}}|f(x_{n})|\\ & =\frac{1}{\sigma_{1}(x_{n})}\left(\sum_{|x|=n}\sigma_{2}^{p}(x)\right)^{{1}/{p}}\geq \frac{1}{\sigma_{1}(x_{n})} \rightarrow\infty
\end{align*}
as $n\rightarrow\infty.$ Therefore, $C_{\varphi}f \notin \mathbb{T}_{\sigma_{1},p}$.
That is, $C_{\varphi}$ is not a bounded operator from $\mathbb{T}_{\sigma_{1},p}$ to $\mathbb{T}_{\sigma_{2},p}$.

Now consider the operator $W_{\psi,\varphi}:\mathbb{T}_{\sigma_{1},p}\rightarrow \mathbb{T}_{\sigma_{2},p}$.
For $f\in \mathbb{T}_{\sigma_{1}, p}$ and $n=2k$, $k\in N_{0}$,
$$
M_{\sigma_{2},p}(n, \psi(f\circ\varphi))=0 \mbox{~as $\psi(x)=0$ if $|x|$ is even.}
$$
For $f\in \mathbb{T}_{\sigma_{1},p}$ and $n=2k-1$, $k\in \mathbb{N},$
\begin{align*}
M_{\sigma_{2},p}(n, \psi(f\circ\varphi)) & =\frac{1}{c_{n}^{ {1}/{p}}}\sigma_{2}(x_{n})|\psi(x_{n})f(\textsl{o})|
 =\frac{|f(\textsl{o})|}{c_{n}^{ {1}/{p}}}\leq |f(\textsl{o})|.
\end{align*}
Therefore, $\|\psi (f \circ \varphi)\|_{\sigma_{2},p}\leq |f(\textsl{o})|
\leq \|f\|$ for all $f\in \mathbb{T}_{\sigma_{1},p}.$
Thus , $W_{\psi, \varphi}:\mathbb{T}_{\sigma_{1},p}\rightarrow \mathbb{T}_
{\sigma_{2},p}$ is a bounded operator.

To see the compactness of $W_{\psi,\varphi}$, let us begin with a bounded sequence
$\{f_{n}\}$ in $\mathbb{T}_{\sigma_{1},p}$ with $f_{n}\rightarrow 0$ pointwise.
As $\|\psi(f_{n}\circ\varphi)\|_{\sigma_{2},p}\leq |f_{n}(\textsl{o})|
\rightarrow 0$ as $n \rightarrow\infty,$
$W_{\psi, \varphi}$ is a compact operator by Theorem \ref{cmpt}.
\eeg

This example shows that $W_{\psi, \varphi}:\mathbb{T}_{\sigma_{1},p}\rightarrow \mathbb{T}_{\sigma_{2},p}$ can be a compact operator even if both the operators $M_{\psi}$ and $C_\varphi$ are unbounded.

As concluding remarks, we leave some unsolved questions:
\begin{itemize}
  \item[(i)] Characterize the bounded, compact, isometric weighted composition
  operators from $\mathbb{T}_{\sigma_{1},p}$ to $\mathbb{T}_{\sigma_{2},p}$, $0<p<\infty$.
  \item[(ii)] The same questions for the operator $W_{\psi, \varphi}:\mathbb{T}_{\sigma_{1},p}\rightarrow \mathbb{T}_{\sigma_{2},q}$,
      $p\neq q$ and $0<p,q<\infty$, are also open.
\end{itemize}

\subsection*{Acknowledgments}
The first author thanks the National Board for Higher Mathematics (NBHM), India, for providing financial support  to carry out this research. The second author is thankful to NBHM/DAE for the project grant,   Grant No. 0211/30/2017/R\&D II/12565. The third author is thankful to CSIR (India), Grant Number 09/1231(0002)/2019-EMR-I.

\subsection*{Conflict of Interests}
The authors declare that there is no conflict of interests regarding the publication of this paper.

\subsection*{Data availability statement}
 The authors declare that this research is purely theoretical and does not associated with any data's.

\end{document}